\documentclass[11pt]{article}
\usepackage{amsmath}
\usepackage{amssymb}
\usepackage{amsfonts}
\usepackage{eucal}
\usepackage{graphicx,url}
\usepackage{amssymb,amsmath, amsthm}
\usepackage{epstopdf}
\def\bincoeff#1#2{{#1\choose #2}}
\newtheorem{theorem}{Theorem}[section]

\newtheorem{example}{Example}[section]

\newtheorem{proposition}[theorem]{Proposition}
\newtheorem{remark}{Remark}

\begin{document}
\title{Correlated fractional counting processes on a finite time interval}
\author{Luisa Beghin\thanks{Address: Dipartimento di
Scienze Statistiche, Sapienza Universit\`a di Roma, Piazzale Aldo
Moro 5, I-00185 Roma, Italy. e-mail:
\texttt{luisa.beghin@uniroma1.it}}\and Roberto
Garra\thanks{Address: Dipartimento di Scienze di Base e Applicate
per l'Ingegneria, Sapienza Universit\`{a} di Roma, Via A. Scarpa
16, I-00161 Roma, Italy. e-mail: \texttt{rolinipame@yahoo.it}}\and
Claudio Macci\thanks{Address: Dipartimento di Matematica,
Universit\`a di Roma Tor Vergata, Via della Ricerca Scientifica,
I-00133 Rome, Italy. e-mail: \texttt{macci@mat.uniroma2.it}}}
\date{}
\maketitle
\begin{abstract}
\noindent We present some correlated fractional counting processes
on a finite time interval. This will be done by considering a
slight generalization of the processes in
\cite{BorgesRodriguesBalakrishnan}. The main case concerns a class
of space-time fractional Poisson processes and, when the
correlation parameter is equal to zero, the univariate
distributions coincide with the ones of the space-time fractional
Poisson process in \cite{OrsingherPolitoSPL2012}. On the other
hand, when we consider the time fractional Poisson process, the
multivariate finite dimensional distributions are different from
the ones presented for the renewal process in
\cite{PolitiKaizojiScalas}. Another case concerns a class of
fractional negative binomial processes.\\

\noindent\emph{Keywords:} Poisson process, negative binomial
process, weighted process.\\
\noindent\emph{Mathematical Subject Classification:} 60G22, 60G55,
60E05, 33E12.
\end{abstract}

\section{Introduction}\label{sec:introduction}
Several fractional processes in the literature are defined by
considering some known equations in terms of suitable fractional
derivatives. In this paper we are interested in particular
L\'{e}vy counting processes, as in the recent paper
\cite{BeghinMacciJAP2014}; in particular we deal with Poisson and
negative binomial processes. There is a wide literature on
fractional Poisson processes: see e.g. \cite{Laskin},
\cite{MainardiGorenfloScalas}, \cite{BeghinOrsingherEJP2009},
\cite{BeghinOrsingherEJP2010}, \cite{OrsingherPolitoSPL2012} and
\cite{PolitiKaizojiScalas} (we also cite
\cite{KumarNaneVellaisamy} and \cite{MeerschaertNaneVellaisamy}
where their representation in terms of randomly time-changed and
subordinated processes was studied in detail). Some references
with fractional negative binomial processes are
\cite{BeghinMacciJAP2014} (see Example 3) and
\cite{VellaisamyMaheshwari}. Among the other fractional processes
in the literature we recall the diffusive processes (see e.g.
\cite{AnguloRuizmedinaAnhGrecksch}, \cite{BaeumerMeerschaertNane},
\cite{Mainardi} \cite{OrsingherBeghinAP2009},
\cite{SchneiderWyss}), the telegraph processes in
\cite{OrsingherBeghinPTRF2004} and the pure birth processes in
\cite{OrsingherPolitoBERNOULLI2010}.

Often the results for these fractional processes are given in
terms of the Mittag-Leffler function
$$E_{\alpha,\beta}(x):=\sum_{r\geq 0}\frac{x^r}{\Gamma(\alpha r+\beta)}$$
(see e.g. \cite{Podlubny}, page 17); we also recalled the
generalized Mittag-Leffler function
$$E_{\alpha,\beta}^\gamma(x):=\sum_{r\geq 0}\frac{(\gamma)^{(r)}x^r}{r!\Gamma(\alpha r+\beta)},$$
where
$$(\gamma)^{(r)}:=\left\{\begin{array}{ll}
\gamma(\gamma+1)\cdots (\gamma+r-1)&\ \mathrm{if}\ r\geq 1\\
1&\ \mathrm{if}\ r=0
\end{array}\right.\ (\mbox{for}\ \gamma\in\mathbb{R})$$
is the rising factorial (also called Pochhammer symbol), and
$E_{\alpha,\beta}^\gamma$ coincides with $E_{\alpha,\beta}$ when
$\gamma=1$.

In this paper we consider some processes
$\{N_\rho(\cdot):\rho\in[0,1]\}$ on a finite time interval
$[0,T]$, for some $T\in (0,\infty)$. More precisely
$N_\rho(\cdot)=\{N_\rho(t):t\in[0,T]\}$ is defined by
$$N_\rho(t):=\sum_{n=1}^{M_g}1_{[0,t]}(X_n^{F,\rho}),$$
where $M_g$ is a nonnegative integer valued random variable with
probability generating function $g$, i.e.
$$g(u):=\mathbb{E}\left[u^{M_g}\right],$$
and $\{X_n^{F,\rho}:n\geq 1\}$ is a sequence of random variables
with (common) distribution function $F$ such that $F(0)=0$ and
$F(T)=1$, and independent of $M_g$; moreover the correlation
coefficient between any pair of random variables $X_n$ and $X_m$,
with $n\neq m$, is equal to a common value $\rho\in[0,1]$.

\begin{remark}\label{rem:t=T}
We have $N_\rho(T)=M_g$; thus the distribution of $N_\rho(T)$ does
not depend on $\rho$.
\end{remark}

In this way we are considering a slight generalization of the
processes presented in \cite{BorgesRodriguesBalakrishnan}; indeed
we can recover several formulas in
\cite{BorgesRodriguesBalakrishnan} by setting
$g(u)=e^{\lambda(u-1)}$ for some $\lambda>0$ (which concerns a
Poisson distributed random variable with mean $\lambda$), and
$F(t)=t$ for $t\in[0,1]$, where $T=1$. The case without
correlation, i.e. the case $\rho=0$, appears in
\cite{BalakrishnanKozubowski}; see also \cite{LefevrePicard} where
that process is considered as a claim number process in insurance.
Here, in view of what follows, we recall the following formulas
(see e.g. (9) and (10) in \cite{BorgesRodriguesBalakrishnan}): we
have the probability generating function
\begin{equation}\label{eq:basic-fgp}
G_{N_\rho(t)}(u)=\rho(1-F(t))+\rho
F(t)g(u)+(1-\rho)g(1-F(t)+F(t)u),
\end{equation}
and the probability mass function
\begin{equation}\label{eq:basic-pmf}
P(N_\rho(t)=k)=(1-\rho)P(N_0(t)=k)+\rho\{(1-F(t))1_{k=0}+F(t)P(M_g=k)\}\
(\mbox{for all}\ k\geq 0),
\end{equation}
where
\begin{equation}\label{eq:basic-pmf-rho=0}
P(N_0(t)=k)=\sum_{n=k}^\infty\bincoeff{n}{k}F^k(t)(1-F(t))^{n-k}P(M_g=n)\
(\mbox{for all}\ k\geq 0)
\end{equation}
concerns the case $\rho=0$ (see (2.4) in
\cite{BalakrishnanKozubowski}).

As pointed out in \cite{BalakrishnanKozubowski}, this class of
counting processes can be useful to tackle the problem of
overdispersion and underdispersion in the analysis of count data
where correlations between events are present. A possible
application can be given for example in models of non-exponential
extinction of radiation in correlated random media (see e.g.
\cite{Kostinski}). We also remark that, as far as the the marginal
distribution of each random variable $N_\rho(t)$, in
\eqref{eq:basic-pmf} we have a mixture between three probability
mass functions, i.e. $\{P(N_0=k):k\geq 0\}$, $\{1_{k=0}:k\geq 0\}$
and $\{P(M_g=k):k\geq 0\}$, and the weights are $1-\rho$,
$\rho(1-F(t))$ and $\rho F(t)$, respectively.

The aim of this paper is to present some correlated fractional
counting processes by choosing in a suitable way the probability
generating function $g$ and a distribution function $F$ above. In
Section \ref{sec:Poisson-processes} we present a class of space-time
fractional Poisson processes (in fact we have the same univariate
distributions of the space-time fractional Poisson process in
\cite{OrsingherPolitoSPL2012} when $\rho=0$). A class of fractional
negative binomial processes is presented in Section
\ref{sec:NB-processes}.

Finally, since the presentation of the results in
\cite{BorgesRodriguesBalakrishnan} refers to the concept of
weighted Poisson processes (see also the previous reference
\cite{BalakrishnanKozubowski} concerning the case $\rho=0$), in
the final Section \ref{sec:weighted-processes} we give some minor
results on weighted processes. Even though this section seems to
be disconnected from the other ones in this paper, in our opinion
it is a small nice enrichment of the content of
\cite{BorgesRodriguesBalakrishnan}.

\section{A class of correlated fractional Poisson processes}\label{sec:Poisson-processes}
For the aims of this section, some preliminaries are needed.
Firstly we consider the Caputo (left fractional) derivative
$\frac{d^\nu}{dt^\nu}$ of order $\nu>0$ (see e.g. ${}^CD_{a+}^\nu$
in (2.4.14) and (2.4.15) in \cite{KilbasSrivastavaTrujillo} with
$a=0$; we use the notation $[x]:=\max\{k\in\mathbb{Z}:k\leq x\}$)
defined by
$$\frac{d^\nu}{dt^\nu}f(t):=\left\{\begin{array}{ll}
\frac{1}{\Gamma(n-\nu)}\int_0^t\frac{1}{(t-s)^{\nu-n+1}}\frac{d^n}{ds^n}f(s)ds&\ \mbox{if}\ \nu\ \mbox{is not integer, where}\ n=[\nu]+1\\
\frac{d^\nu}{dt^\nu}f(t)&\ \mbox{if}\ \nu\ \mbox{is integer}
\end{array}\right.\ (\mbox{for all}\ t\geq 0);$$
note that, since here we consider $\nu\in(0,1]$, we have (see e.g.
(2.4.17) in \cite{KilbasSrivastavaTrujillo} with $a=0$)
$$\frac{d^\nu}{dt^\nu}f(t):=\left\{\begin{array}{ll}
\frac{1}{\Gamma(1-\nu)}\int_0^t\frac{1}{(t-s)^{\nu}}\frac{d}{ds}f(s)ds&\ \mbox{if}\ \nu\in(0,1)\\
\frac{d}{dt}f(t)&\ \mbox{if}\ \nu=1
\end{array}\right.\ (\mbox{for all}\ t\geq 0).$$
We also consider the (fractional) difference operator $(I-B)^\alpha$
in \cite{OrsingherPolitoSPL2012}; more precisely $I$ is the identity
operator, $B$ is the backward shift operator defined by
$Bf(k)=f(k-1)$ and $B^{r-1}Bf(k)=f(k-r)$, and therefore
\begin{equation}\label{eq:fractional-version-of-I-B}
(I-B)^\alpha=\sum_{j=0}^\infty(-1)^j\bincoeff{\alpha}{j}B^j.
\end{equation}
We now recall that Orsingher and Polito in
\cite{OrsingherPolitoSPL2012} considered the space-time fractional
Poisson process $\{N_0^{\alpha,\nu}(t):t\geq 0\}$, for
$\alpha,\nu\in(0,1]$, whose probability mass functions
$\{p_k(t):k\geq 0\}$ solve the Cauchy problem
\begin{equation}\nonumber
\begin{cases}
\frac{d^{\nu}} {dt^{\nu}}p_k(t)= -\lambda^{\alpha}
(I-B)^{\alpha}p_k(t)\\
p_k(0)= \begin{cases} 0, \quad k>0, \\
1, \quad k= 0.
\end{cases}
\end{cases}
\end{equation}
The explicit form of the probability generating function of this
process has the form (see \cite{OrsingherPolitoSPL2012}, equation
(2.28))
$$\mathbb{E}\left[u^{N_0^{\alpha,\nu}(t)} \right]=E_{\nu,1}(-\lambda^{\alpha} t^{\nu}(1-u)^{\alpha}).$$

In this section we consider a class of correlated space-time
fractional Poisson processes on a finite time interval $[0,T]$.
For $\alpha,\nu\in(0,1]$ we consider
$N_\rho(\cdot)=N_\rho^{\alpha,\nu}(\cdot)$ such that the
probability generating function of $M_g$ is
\begin{equation}\label{eq:spacetimefractionalPoisson-fgp}
g(u):=E_{\nu,1}(-\lambda^\alpha T^\nu(1-u)^\alpha),
\end{equation}
and the distribution function of the random variables
$\{X_n^{F,\rho}:n\geq 1\}$ is
$$F(t):=(t/T)^{\nu/\alpha}\ (\mbox{for}\ t\in [0,T]).$$

In what follows we present the probability generating functions in
Proposition \ref{prop:spacetimefractionalPoisson-pgf} and the
corresponding probability mass functions in Proposition
\ref{prop:spacetimefractionalPoisson-pmf}. Moreover, in
Proposition \ref{prop:spacetimefractionalPoisson-eq-for-pmf}, we
give an equation for the probability mass functions in Proposition
\ref{prop:spacetimefractionalPoisson-pmf} with respect to time
$t$.

\begin{proposition}\label{prop:spacetimefractionalPoisson-pgf}
The probability generating functions
$\{G_{N_\rho^{\alpha,\nu}(t)}:t\in [0,T]\}$ are
\begin{align*}
G_{N_\rho^{\alpha,\nu}(t)}(u)&=\rho(1-(t/T)^{\nu/\alpha})+\rho(t/T)^{\nu/\alpha}E_{\nu,1}(-\lambda^\alpha T^\nu(1-u)^\alpha)\\
&+(1-\rho)E_{\nu,1}(-\lambda^\alpha t^\nu(1-u)^\alpha).
\end{align*}
\end{proposition}
\noindent\emph{Proof.} We have
\begin{align*}
G_{N_\rho^{\alpha,\nu}(t)}(u)&=\rho(1-(t/T)^{\nu/\alpha})+\rho(t/T)^{\nu/\alpha}E_{\nu,1}(-\lambda^\alpha T^\nu(1-u)^\alpha)\\
&+(1-\rho)E_{\nu,1}(-\lambda^\alpha
T^\nu(1-\{1-(t/T)^{\nu/\alpha}+(t/T)^{\nu/\alpha}u\})^\alpha)
\end{align*}
by \eqref{eq:basic-fgp}, and we conclude with some manipulations
of the last term. $\Box$

\begin{remark}\label{rem:comparison-with-renewal-case}
By Proposition \ref{prop:spacetimefractionalPoisson-pgf}, if
$\rho=0$ we have the probability generating function
\begin{equation}\label{eq:spacetimefractionalPoisson-pgf-rho=0}
G_{N_0^{\alpha,\nu}(t)}(u)=E_{\nu,1}(-\lambda^\alpha
t^\nu(1-u)^\alpha)
\end{equation}
which coincides with the one presented in the last case of Table 1
in \cite{OrsingherPolitoSPL2012}; note that
\eqref{eq:spacetimefractionalPoisson-pgf-rho=0} is a generalization
of \eqref{eq:spacetimefractionalPoisson-fgp} with $t\in[0,T]$
instead of $t=T$. Thus the univariate distributions of the random
variables $\{N_0^{1,\nu}(t):t\in[0,T]\}$ (for the case $\alpha=1$)
coincide with the ones of the random variables of the renewal
process $\{M(t):t\in[0,T]\}$ in \cite{PolitiKaizojiScalas}
(restricted to the same finite time interval). On the other hand one
can check that the multivariate finite dimensional marginal
distributions are different from the ones in
\cite{PolitiKaizojiScalas} (and, in fact,
$\{N_\rho^{\alpha,\nu}(t):t\in[0,T]\}$ is not a renewal process). We
explain this with a simple example where we take into account that
$$P(M(s)=1)=P(N_0^{1,\nu}(s)=1)=\lambda s^\nu E_{\nu,\nu+1}^2(-\lambda s^\nu)\ (\mbox{for}\ s\in[0,T])$$
by (2.5) in \cite{BeghinOrsingherEJP2010}. In fact, for
$t\in(0,T)$, we have
\begin{equation}\label{eq:conditional-probability-KPS}
P(M(t)=1,M(T)=1)=\lambda t^\nu E_{\nu,\nu+1}^2(-\lambda
t^\nu)E_{\nu,1}(-\lambda (T-t)^\nu)
\end{equation}
by combining (11) and (14) in \cite{PolitiKaizojiScalas} (with
$(t_1,t_2)=(t,T)$ and $(n_1,n_2)=(1,1)$) with (2) and (4) in the
same reference, and
\begin{equation}\label{eq:conditional-probability-BRB}
P(N_0^{1,\nu}(t)=1,N_0^{1,\nu}(T)=1)=\frac{t}{T}\lambda T^\nu
E_{\nu,\nu+1}^2(-\lambda T^\nu)
\end{equation}
because
$P(N_0^{\alpha,\nu}(t)=1|N_0^{\alpha,\nu}(T)=1)=\frac{t}{T}$ by
construction. Then \eqref{eq:conditional-probability-KPS} and
\eqref{eq:conditional-probability-BRB} coincide only for the
non-fractional case $\nu=1$ (see Figure \ref{fig} below).
\end{remark}
\begin{figure}[ht]
\begin{center}
\includegraphics[angle=0,width=0.5\textwidth]{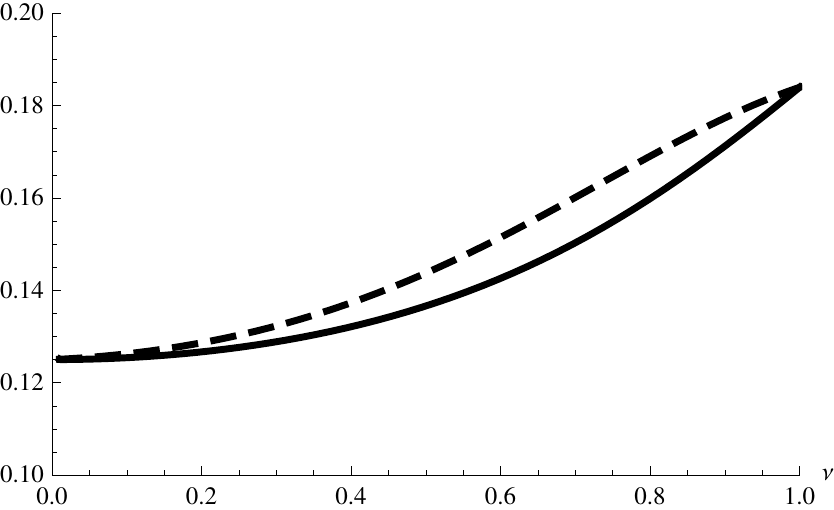}
\caption{The probabilities \eqref{eq:conditional-probability-KPS}
(dashed line) and \eqref{eq:conditional-probability-BRB} (solid
line) versus $\nu\in(0,1]$ for $t=1/2$ and
$T=\lambda=1$.}\label{fig}
\end{center}
\end{figure}

\begin{proposition}\label{prop:spacetimefractionalPoisson-pmf}
The probability mass functions
$\{P(N_\rho^{\alpha,\nu}(t)=\cdot):t\in [0,T]\}$ are
\begin{align*}
P(N_\rho^{\alpha,\nu}(t)=k)&=(1-\rho)\frac{(-1)^k}{k!}\sum_{r=0}^\infty\frac{(-\lambda^\alpha
t^\nu)^r}{\Gamma(\nu r+1)}\frac{\Gamma(\alpha r+1)}{\Gamma(\alpha
r+1-k)}+\rho\left\{\left(1-(t/T)^{\frac{\nu}{\alpha}}\right)1_{k=0}\right.\\
&\left.+(t/T)^{\frac{\nu}{\alpha}}\cdot\frac{(-1)^k}{k!}\sum_{r=0}^\infty\frac{(-\lambda^\alpha
T^\nu)^r}{\Gamma(\nu r+1)}\frac{\Gamma(\alpha r+1)}{\Gamma(\alpha
r+1-k)}\right\}\ (\mbox{for all}\ k\geq 0).
\end{align*}
\end{proposition}
\noindent\emph{Proof.} Firstly we have
\begin{equation}\label{eq:spacetimefractionalPoisson-pmf}
P(M_g=n)=P(N_\rho^{\alpha,\nu}(T)=k)=\frac{(-1)^n}{n!}\sum_{r=0}^\infty\frac{(-\lambda^\alpha
T^\nu)^r}{\Gamma(\nu r+1)}\frac{\Gamma(\alpha r+1)}{\Gamma(\alpha
r+1-n)}\ (\mbox{for all}\ n\geq 0)
\end{equation}
by the probability generating function in
\eqref{eq:spacetimefractionalPoisson-fgp} (see (1.8) in
\cite{OrsingherPolitoSPL2012}) and by Remark \ref{rem:t=T}.
Moreover, if we consider \eqref{eq:basic-pmf-rho=0}, we get
\begin{align*}
P(N_0^{\alpha,\nu}(t)=k)&=\sum_{n=k}^\infty\bincoeff{n}{k}(t/T)^{\frac{\nu}{\alpha}k}
\left(1-(t/T)^{\frac{\nu}{\alpha}}\right)^{n-k}\frac{(-1)^n}{n!}\sum_{r=0}^\infty\frac{(-\lambda^\alpha
T^\nu)^r}{\Gamma(\nu r+1)}\frac{\Gamma(\alpha r+1)}{\Gamma(\alpha
r+1-n)}\\
&=\frac{(-1)^k}{k!}(t/T)^{\frac{\nu}{\alpha}k}\sum_{n=k}^\infty\frac{(-1)^{n-k}}{(n-k)!}
\left(1-(t/T)^{\frac{\nu}{\alpha}}\right)^{n-k}\\
&\cdot\sum_{r=0}^\infty\frac{(-\lambda^\alpha t^\nu)^r}{\Gamma(\nu
r+1)}(T/t)^{\nu r}\frac{\Gamma(\alpha r+1)}{\Gamma(\alpha
r+1-k)}\frac{\Gamma(\alpha
r+1-k)}{\Gamma(\alpha r+1-n)}\\
&=\frac{(-1)^k}{k!}(t/T)^{\frac{\nu}{\alpha}k}\sum_{r=0}^\infty\frac{(-\lambda^\alpha
t^\nu)^r}{\Gamma(\nu r+1)}(T/t)^{\nu
r}\frac{\Gamma(\alpha r+1)}{\Gamma(\alpha r+1-k)}\\
&\cdot\sum_{j=0}^\infty\frac{(-1)^j}{j!}
\left(1-(t/T)^{\frac{\nu}{\alpha}}\right)^j\frac{\Gamma(\alpha
r+1-k)}{\Gamma(\alpha r+1-k-j)}\ (\mbox{for all}\ k\geq 0).
\end{align*}
Then, by a well-known \lq\lq Newton's generalized binomial
theorem\rq\rq, we obtain
\begin{align*}
P(N_0^{\alpha,\nu}(t)=k)&=\frac{(-1)^k}{k!}(t/T)^{\frac{\nu}{\alpha}k}\sum_{r=0}^\infty\frac{(-\lambda^\alpha
t^\nu)^r}{\Gamma(\nu r+1)}(T/t)^{\nu r}\frac{\Gamma(\alpha
r+1)}{\Gamma(\alpha
r+1-k)}\left(1+(t/T)^{\frac{\nu}{\alpha}}-1\right)^{\alpha
r-k}\\
&=\frac{(-1)^k}{k!}\sum_{r=0}^\infty\frac{(-\lambda^\alpha
t^\nu)^r}{\Gamma(\nu r+1)}(t/T)^{\frac{\nu}{\alpha}k-\nu r+\nu
r-\frac{\nu}{\alpha}k}\frac{\Gamma(\alpha r+1)}{\Gamma(\alpha
r+1-k)}\\
&=\frac{(-1)^k}{k!}\sum_{r=0}^\infty\frac{(-\lambda^\alpha
t^\nu)^r}{\Gamma(\nu r+1)}\frac{\Gamma(\alpha r+1)}{\Gamma(\alpha
r+1-k)}\ (\mbox{for all}\ k\geq 0)
\end{align*}
where, as we expected by
\eqref{eq:spacetimefractionalPoisson-pgf-rho=0},
$P(N_0^{\alpha,\nu}(t)=k)$ here meets
$P(N_\rho^{\alpha,\nu}(T)=k)$ in
\eqref{eq:spacetimefractionalPoisson-pmf} (here we have $t$ and
$k$ in place of $T$ and $n$ in
\eqref{eq:spacetimefractionalPoisson-pmf}). We conclude the proof
by considering \eqref{eq:basic-pmf} together with the last above
expression obtained for the case $\rho=0$. $\Box$\\

In view of the next Proposition
\ref{prop:spacetimefractionalPoisson-eq-for-pmf} we remark that in a
part of the proof we refer to Theorem 2 in \cite{BeghinDovidio}
which can be derived by referring to a subordinated representation
of the space-time fractional Poisson process in terms of both stable
subordinator and its inverse (see also (3.20), together with (3.1),
in the same reference).

\begin{proposition}\label{prop:spacetimefractionalPoisson-eq-for-pmf}
Let $\{P(N_\rho^{\alpha,\nu}(t)=\cdot):t\in [0,T]\}$ be the
probability mass functions in Proposition
\ref{prop:spacetimefractionalPoisson-pmf}. Then we have the
following equations: for $k=0$,
\begin{align*}
\frac{d^\nu}{dt^\nu}P(N_\rho^{\alpha,\nu}(t)=0)&=-\lambda^\alpha
P(N_\rho^{\alpha,\nu}(t)=0)+\lambda^\alpha\rho\\
&-\rho(t/T)^{\nu/\alpha}\left(\lambda^\alpha+t^{-\nu}\frac{\Gamma(\nu/\alpha+1)}{\Gamma(\nu/\alpha-\nu+1)}\right)(1-P(N_\rho^{\alpha,\nu}(T)=0));
\end{align*}
for all $k\geq 1$,
\begin{align*}
\frac{d^\nu}{dt^\nu}P(N_\rho^{\alpha,\nu}(t)=k)=&-\lambda^\alpha(I-B)^\alpha
P(N_\rho^{\alpha,\nu}(t)=k)+\lambda^\alpha\rho(1-(t/T)^{\nu/\alpha})(-1)^k\bincoeff{\alpha}{k}\\
&+\rho(t/T)^{\nu/\alpha}\left[\lambda^\alpha(I-B)^\alpha+t^{-\nu}\frac{\Gamma(\nu/\alpha+1)}{\Gamma(\nu/\alpha-\nu+1)}\right]P(N_\rho^{\alpha,\nu}(T)=k).
\end{align*}
In all cases we have the initial conditions
$P(N_\rho^{\alpha,\nu}(0)=0)=1$ and
$P(N_\rho^{\alpha,\nu}(0)=k)=0$ for all $k\geq 1$.
\end{proposition}
\noindent\emph{Proof.} The initial conditions trivially hold.
Throughout this proof we consider the notation
$$p_k^{\alpha,\nu}(t)=P(N_0^{\alpha,\nu}(t)=k)\ (\mbox{for all}\ k\geq 0)$$
for the probability mass function concerning the case $\rho=0$.
Then, by \eqref{eq:basic-pmf} and Remark \ref{rem:t=T}, we get
\begin{align*}
\frac{d^\nu}{dt^\nu}P(N_\rho^{\alpha,\nu}(t)=k)=&(1-\rho)\frac{d^\nu}{dt^\nu}p_k^{\alpha,\nu}(t)
+\rho\left\{-\frac{1}{T^{\nu/\alpha}}1_{k=0}\frac{d^\nu}{dt^\nu}t^{\nu/\alpha}
+\frac{1}{T^{\nu/\alpha}}P(N_\rho^{\alpha,\nu}(T)=k)\frac{d^\nu}{dt^\nu}t^{\nu/\alpha}\right\}\\
=&(1-\rho)\frac{d^\nu}{dt^\nu}p_k^{\alpha,\nu}(t)
-\frac{\rho}{T^{\nu/\alpha}}\left\{1_{k=0}-P(N_\rho^{\alpha,\nu}(T)=k)\right\}\frac{d^\nu}{dt^\nu}t^{\nu/\alpha}.
\end{align*}
Moreover we have
$$\frac{d^\nu}{dt^\nu}p_k^{\alpha,\nu}(t)=-\lambda^\alpha(I-B)^\alpha p_k^{\alpha,\nu}(t)$$
by Theorem 2 in \cite{BeghinDovidio} and
$$\frac{d^\nu}{dt^\nu}t^{\nu/\alpha}=t^{\nu/\alpha-\nu}\frac{\Gamma(\nu/\alpha+1)}{\Gamma(\nu/\alpha-\nu+1)}$$
(see e.g. (2.2.11) and (2.4.8) in \cite{KilbasSrivastavaTrujillo},
or a correction of (2.4.28) in the same reference); then we get
\begin{align*}
\frac{d^\nu}{dt^\nu}P(N_\rho^{\alpha,\nu}(t)=k)&=-\lambda^\alpha(I-B)^\alpha(1-\rho)p_k^{\alpha,\nu}(t)\\
&-\rho(t/T)^{\nu/\alpha}\left\{1_{k=0}-P(N_\rho^{\alpha,\nu}(T)=k)\right\}t^{-\nu}\frac{\Gamma(\nu/\alpha+1)}{\Gamma(\nu/\alpha-\nu+1)}.
\end{align*}
From now on we consider the cases $k=0$ and $k\geq 1$
separately.\\
\emph{Case} $k=0$. Firstly we have
$$(I-B)^\alpha p_0^{\alpha,\nu}(t)=\sum_{j=0}^\infty(-1)^j\bincoeff{\alpha}{j}p_{0-j}^{\alpha,\nu}(t)=p_0^{\alpha,\nu}(t)$$
by \eqref{eq:fractional-version-of-I-B}; therefore
\begin{align*}
\frac{d^\nu}{dt^\nu}P(N_\rho^{\alpha,\nu}(t)=0)=-\lambda^\alpha(1-\rho)p_0^{\alpha,\nu}(t)
-\rho(t/T)^{\nu/\alpha}\left\{1-P(N_\rho^{\alpha,\nu}(T)=0)\right\}t^{-\nu}\frac{\Gamma(\nu/\alpha+1)}{\Gamma(\nu/\alpha-\nu+1)}.
\end{align*}
then, by \eqref{eq:basic-pmf} and Remark \ref{rem:t=T},
\begin{align*}
\frac{d^\nu}{dt^\nu}P(N_\rho^{\alpha,\nu}(t)=0)&=-\lambda^\alpha\{P(N_\rho^{\alpha,\nu}(t)=0)
-\rho\{1-(t/T)^{\nu/\alpha}+(t/T)^{\nu/\alpha}P(N_\rho^{\alpha,\nu}(T)=0)\}\}\\
&-\rho(t/T)^{\nu/\alpha}\left\{1-P(N_\rho^{\alpha,\nu}(T)=0)\right\}t^{-\nu}\frac{\Gamma(\nu/\alpha+1)}{\Gamma(\nu/\alpha-\nu+1)}.
\end{align*}
and, finally, we can check by inspection that the last equation is
equivalent to the one in the statement of the proposition.\\
\emph{Case} $k\geq 1$. Firstly, again by \eqref{eq:basic-pmf} and
Remark \ref{rem:t=T}, we have
\begin{align*}
\frac{d^\nu}{dt^\nu}P(N_\rho^{\alpha,\nu}(t)=k)&=-\lambda^\alpha(I-B)^\alpha[P(N_\rho^{\alpha,\nu}(t)=k)\\
&-\rho(1-(t/T)^{\nu/\alpha})1_{k=0}-\rho(t/T)^{\nu/\alpha}P(N_\rho^{\alpha,\nu}(T)=k)]\\
&+\rho(t/T)^{\nu/\alpha}P(N_\rho^{\alpha,\nu}(T)=k)t^{-\nu}\frac{\Gamma(\nu/\alpha+1)}{\Gamma(\nu/\alpha-\nu+1)}\\
&=-\lambda^\alpha(I-B)^\alpha
P(N_\rho^{\alpha,\nu}(t)=k)+\lambda^\alpha\rho(1-(t/T)^{\nu/\alpha})(I-B)^\alpha1_{k=0}\\
&+\rho(t/T)^{\nu/\alpha}\left[\lambda^\alpha(I-B)^\alpha+t^{-\nu}\frac{\Gamma(\nu/\alpha+1)}{\Gamma(\nu/\alpha-\nu+1)}\right]P(N_\rho^{\alpha,\nu}(T)=k);
\end{align*}
then we get the desired equation by noting that
$$(I-B)^\alpha 1_{k=0}=\sum_{j=0}^\infty(-1)^j\bincoeff{\alpha}{j}1_{k-j=0}=(-1)^k\bincoeff{\alpha}{k}.$$
The proof is complete. $\Box$\\

Finally we remark that, even if the equations in Proposition
\ref{prop:spacetimefractionalPoisson-eq-for-pmf} have some
analogies with other results for fractional Poisson processes in
the literature, here some standard techniques do not work because
we deal with a finite horizon time case (i.e. $t\in [0,T]$).

\section{A class of correlated fractional negative binomial processes}\label{sec:NB-processes}
It is well-known that the negative binomial process can be seen as
a suitable compound Poisson process with logarithmic distributed
summands (see e.g. Proposition 1.1 in \cite{KozubowskiPodgorski}).
More precisely, for some $p\in(0,1)$ and some integer $r\geq 1$,
we have the probability generating function
$$u\mapsto h^r(m(u)),$$
where:
$$h(u):=e^{\lambda(u-1)},\ \mbox{with}\ \lambda=-\log p,$$
is the probability generating function of a Poisson distributed
random variable with mean $\lambda=-\log p$;
$$m(u):=\frac{\log(1-(1-p)u)}{\log p},\ \mbox{for}\ |u|<\frac{1}{1-p},$$
is the probability generating function of a logarithmic
distributed random variable (obviously we have $m(u)=\infty$ if
$|u|\geq\frac{1}{1-p}$).

In this section we present a class of correlated fractional
negative binomial processes on a finite time interval $[0,T]$.
More precisely we consider the same approach with the probability
generating function of a space-time fractional Poisson distributed
random variable; thus, for $\alpha,\nu\in(0,1]$, we have
$$h_{\alpha,\nu}(u):=E_{\nu,1}(-\lambda^\alpha(1-u)^\alpha)$$
in place of $h$ (note that $h$ coincides with $h_{1,1}$), again
with $\lambda=-\log p$, and this meets $g$ in
\eqref{eq:spacetimefractionalPoisson-fgp} with $T=1$. Thus we have
\begin{equation}\label{eq:spacetimefractionalNB-pgf}
g(u):=\left\{E_{\nu,1}\left(-(-\log
p)^\alpha\left(1-\frac{\log(1-(1-p)u)}{\log
p}\right)^\alpha\right)\right\}^r
=\left\{E_{\nu,1}\left(-\log^\alpha\left(\frac{1-(1-p)u}{p}\right)\right)\right\}^r,
\end{equation}
where, again, $r\geq 1$ is an integer power of the function
$E_{\nu,1}$, $p\in(0,1)$ and $|u|<\frac{1}{1-p}$. We remark that
$g$ in \eqref{eq:spacetimefractionalNB-pgf} is the probability
generating function of $N_\rho(T)$, but it does not depend on $T$
as happens for $g$ in \eqref{eq:spacetimefractionalPoisson-fgp}.

As far as the distribution function $F$ is concerned, we argue as
in Section \ref{sec:Poisson-processes} as follows: for all
$t\in[0,T]$, we want to have the condition
$$G_{N_0^{\alpha,\nu}(t)}(u)=\left\{E_{\nu,1}\left(-\log^\alpha\left(\frac{1-(1-q(t))u}{q(t)}\right)\right)\right\}^r$$
for some $q(\cdot)$ such that $q(t)\in(0,1]$ for all $t\in[0,T]$
and $q(T)=p$. Then, by \eqref{eq:basic-fgp} with $\rho=0$ and by
\eqref{eq:spacetimefractionalNB-pgf}, we require that
\begin{align*}
\frac{1-(1-q(t))u}{q(t)}=&\frac{1-(1-p)(1-F(t)+F(t)u)}{p}\\
=&\frac{1-(1-p)(1-F(t))-(1-p)F(t)u}{p};
\end{align*}
so, if we divide both numerator and denominator by
$1-(1-p)(1-F(t))$, we get
$$q(t)=\frac{p}{1-(1-p)(1-F(t))}.$$
Moreover we have
$$q(t)=\frac{p}{p+(1-p)F(t)}=\frac{1}{1+(\frac{1}{p}-1)F(t)}$$
which yields
\begin{equation}\label{eq:dfNB}
F(t):=\frac{\frac{1}{q(t)}-1}{\frac{1}{p}-1}\ (\mbox{for}\ t\in
[0,T]),
\end{equation}
and the function $q(\cdot)$ has to be a decreasing. We also give a
particular example with a choice of $q(\cdot)$, and we provide the
corresponding distribution function $F$.

\begin{example}\label{ex:NB-ex-BK}
If we set
$$q(t)=\frac{1-\lambda}{1-(1-\frac{t}{T})\lambda}$$
for some $\lambda\in(0,1)$, we recover the example in Section 3.3
in \cite{BalakrishnanKozubowski} (see also Section 4.3 in
\cite{BorgesRodriguesBalakrishnan} for a generalization). In fact
this choice of $q(\cdot)$ is the analogue of (3.6) in
\cite{BalakrishnanKozubowski}; moreover, if we set $p=1-\lambda$,
we have
$$q(t)=\frac{p}{1-(1-\frac{t}{T})(1-p)}=\frac{1}{1+(\frac{1}{p}-1)\frac{t}{T}}$$
and therefore $F(t)=\frac{t}{T}$.
\end{example}

In what follows we present the probability generating functions in
Proposition \ref{prop:spacetimefractionalNB-pgf} and, for $r=1$
only, the corresponding probability mass functions in Proposition
\ref{prop:spacetimefractionalNB-pmf} (for $r\geq 2$ we have the
$r$-th convolution of the probability mass function of the case
$r=1$, but we cannot provide manageable formulas). Moreover, in
Proposition \ref{prop:spacetimefractionalNB-eq-for-pgf}, we give
an equation for the probability generating functions
$\{G_{N_\rho^{\alpha,\nu}(t)}:t\in [0,T]\}$ in Proposition
\ref{prop:spacetimefractionalNB-pgf} for $r=1$, $\nu=\alpha$ and
$\rho\in\{0,1\}$; in this case we consider fractional derivatives
with respect to their argument $u$, and not with respect to time
$t$.

\begin{proposition}\label{prop:spacetimefractionalNB-pgf}
The probability generating functions
$\{G_{N_\rho^{\alpha,\nu}(t)}:t\in [0,T]\}$ are
\begin{align*}
G_{N_\rho^{\alpha,\nu}(t)}(u)&=\rho\left(1-\frac{\frac{1}{q(t)}-1}{\frac{1}{p}-1}\right)+\rho
\frac{\frac{1}{q(t)}-1}{\frac{1}{p}-1}\left\{E_{\nu,1}\left(-\log^\alpha\left(\frac{1-(1-p)u}{p}\right)\right)\right\}^r\\
&+(1-\rho)\left\{E_{\nu,1}\left(-\log^\alpha\left(\frac{1-(1-q(t))u}{q(t)}\right)\right)\right\}^r.
\end{align*}
\end{proposition}
\noindent\emph{Proof.} This is an immediate consequence of
\eqref{eq:basic-fgp} and the formulas above. $\Box$\\

In view of the next Proposition
\ref{prop:spacetimefractionalNB-pmf} some preliminaries are
needed. Firstly we consider the Stirling numbers $\{s_{k,h}:k\geq
h\geq 0\}$; for their definition and some properties used below
see e.g. \cite{AbramowitzStegun}, page 824. Moreover
$$_p\Psi_q\left[\begin{array}{cc}
(a_1,\alpha_1)\ldots(a_p,\alpha_p)\\
(b_1,\beta_1)\ldots(b_q,\beta_q)
\end{array}\right](z):=\sum_{j\geq 0}\frac{\prod_{h=1}^p\Gamma(a_h+\alpha_h j)}{\prod_{k=1}^q\Gamma(b_k+\beta_k j)}\frac{z^j}{j!}$$
is the Fox-Wright function (see e.g. (1.11.14) in
\cite{KilbasSrivastavaTrujillo}) under the convergence condition
\begin{equation}\label{eq:convergence-condition-wright}
\sum_{k=1}^q\beta_k-\sum_{h=1}^p\alpha_h>-1
\end{equation}
(see e.g. (1.11.15) in \cite{KilbasSrivastavaTrujillo}).

\begin{proposition}\label{prop:spacetimefractionalNB-pmf}
If $r=1$, the probability mass functions
$\{P(N_\rho^{\alpha,\nu}(t)=\cdot):t\in [0,T]\}$ are
\begin{align*}
P(N_\rho^{\alpha,\nu}(t)=k)&=(1-\rho)P(N_0^{\alpha,\nu}(t)=k)\\
&+\rho\left\{\frac{\frac{1}{p}-\frac{1}{q(t)}}{\frac{1}{p}-1}1_{k=0}
+\frac{\frac{1}{q(t)}-1}{\frac{1}{p}-1}P(N_0^{\alpha,\nu}(T)=k)\right\}\
(\mbox{for all}\ k\geq 0),
\end{align*}
where, for all $t\in[0,T]$,
\begin{equation}\label{eq:spacetimefractionalNB-pgf-rho=0}
P(N_0^{\alpha,\nu}(t)=k)=\left\{\begin{array}{ll}
E_{\nu,1}\left(-\log^\alpha\left(1+A_t\right)\right)&\ \mbox{if}\ k=0\\
\frac{1}{k!}\frac{(-A_t)^k}{(1+A_t)^k}\sum_{h=1}^k\log^{-h}\left(1+A_t\right)s_{k,h}&\ \\
\cdot\ _2\Psi_2\left[\begin{array}{cc}
(1,\alpha)&(1,1)\\
(1-h,\alpha)&(1,\nu)
\end{array}\right](-\log^\alpha(1+A_t))&\ \mbox{if}\ k\geq 1
\end{array}\right.
\end{equation}
and $A_t:=\frac{1}{q(t)}-1$ (note that the convergence condition
\eqref{eq:convergence-condition-wright} holds because we have
$\alpha+\nu-(\alpha+1)=\nu-1>-1$).
\end{proposition}
\noindent\emph{Proof.} Firstly we remark that we can only check
\eqref{eq:spacetimefractionalNB-pgf-rho=0} (concerning the case
$\rho=0$); in fact we obtain the formula for the general case by
combining \eqref{eq:basic-pmf}, $F$ in \eqref{eq:dfNB} and
\eqref{eq:spacetimefractionalNB-pgf-rho=0}. It is well-known that
\begin{equation}\label{eq:well-known-formula}
P(N_0^{\alpha,\nu}(t)=k)=\left\{\begin{array}{ll}
G_{N_0^{\alpha,\nu}(t)}(0)&\ \mbox{if}\ k=0\\
\frac{1}{k!}\left.\frac{d^k}{du^k}G_{N_0^{\alpha,\nu}(t)}(u)\right|_{u=0}&\
\mbox{if}\ k\geq 1.
\end{array}\right.
\end{equation}
Firstly, if $A_t=\frac{1}{q(t)}-1$ as in the statement of the
proposition, we have
$$G_{N_0^{\alpha,\nu}(t)}(u)=E_{\nu,1}\left(-\log^\alpha\left(\frac{1-(1-q(t))u}{q(t)}\right)\right)=E_{\nu,1}\left(-\log^\alpha\left(1+A_t(1-u)\right)\right),$$
and we immediately obtain
\eqref{eq:spacetimefractionalNB-pgf-rho=0} for $k=0$. Moreover, if
we prove that
\begin{align}
\nonumber\frac{d^k}{du^k}&E_{\nu,1}\left(-\log^\alpha\left(1+A(1-u)\right)\right)\\
\label{eq:derivative-by-induction}&=\frac{(-A)^k}{(1+A(1-u))^k}\sum_{j\geq
0}\sum_{h=1}^k\frac{(-1)^j\Gamma(\alpha j+1)}{\Gamma(\alpha
j-h+1)\Gamma(\nu j+1)}s_{k,h}\log^{\alpha
j-h}\left(1+A(1-u)\right)\\
\label{eq:derivative-with-wright-function}&=\frac{1}{k!}\frac{(-A)^k}{(1+A)^k}\sum_{h=1}^k\log^{-h}\left(1+A\right)s_{k,h}\cdot\
_2\Psi_2\left[\begin{array}{cc}
(1,\alpha)&(1,1)\\
(1-h,\alpha)&(1,\nu)
\end{array}\right](-\log^\alpha(1+A));
\end{align}
for $k\geq 1$ (and for all $A\in\mathbb{R}$), we obtain
\eqref{eq:spacetimefractionalNB-pgf-rho=0} for $k\geq 1$ (and the
proof is complete) as an immediate consequence of
\eqref{eq:well-known-formula} and
\eqref{eq:derivative-with-wright-function} with $A=A_t$. Therefore
in the remaining part of the proof we only prove the first
equality \eqref{eq:derivative-by-induction} by induction; in fact
the second equality \eqref{eq:derivative-with-wright-function} can
be checked by inspection. For $k=1$ we have
$$\frac{d}{du}E_{\nu,1}\left(-\log^\alpha\left(1+A(1-u)\right)\right)=\sum_{j\geq 0}\frac{(-1)^j\alpha j}{\Gamma(\nu
j+1)}\frac{\log^{\alpha
j-1}\left(1+A(1-u)\right)}{1+A(1-u)}\cdot(-A),$$ and
\eqref{eq:derivative-by-induction} is proved noting that
$s_{1,1}=1$ and $\alpha j=\frac{\Gamma(\alpha j+1)}{\Gamma(\alpha
j)}$. Now we assume that \eqref{eq:derivative-by-induction} holds
for $k>1$. Then we have
\begin{align*}
\frac{d^{k+1}}{du^{k+1}}&E_{\nu,1}\left(-\log^\alpha\left(1+A(1-u)\right)\right)\\
&=\frac{d}{du}\left\{\frac{(-A)^k}{(1+A(1-u))^k}\sum_{j\geq
0}\sum_{h=1}^k\frac{(-1)^j\Gamma(\alpha j+1)}{\Gamma(\alpha
j-h+1)\Gamma(\nu j+1)}s_{k,h}\log^{\alpha
j-h}\left(1+A(1-u)\right)\right\}\\
&=(-A)^k\left\{\frac{(-k)(-A)}{(1+A(1-u))^{k+1}}\sum_{j\geq
0}\sum_{h=1}^k\frac{(-1)^j\Gamma(\alpha j+1)}{\Gamma(\alpha
j-h+1)\Gamma(\nu j+1)}s_{k,h}\log^{\alpha
j-h}\left(1+A(1-u)\right)\right.\\
&+\left.\frac{1}{(1+A(1-u))^k}\sum_{j\geq
0}\sum_{h=1}^k\frac{(-1)^j\Gamma(\alpha j+1)}{\Gamma(\alpha
j-h+1)\Gamma(\nu j+1)}s_{k,h}\frac{(\alpha j-h)\log^{\alpha
j-h-1}\left(1+A(1-u)\right)}{1+A(1-u)}\cdot(-A)\right\},
\end{align*}
and we obtain
\begin{align*}
\frac{d^{k+1}}{du^{k+1}}&E_{\nu,1}\left(-\log^\alpha\left(1+A(1-u)\right)\right)\\
=&\frac{(-A)^{k+1}}{(1+A(1-u))^{k+1}}\left\{-k\sum_{j\geq
0}\sum_{h=1}^k\frac{(-1)^j\Gamma(\alpha j+1)}{\Gamma(\alpha
j-h+1)\Gamma(\nu j+1)}s_{k,h}\log^{\alpha
j-h}\left(1+A(1-u)\right)\right.\\
&+\left.\sum_{j\geq 0}\sum_{h=0}^k\frac{(-1)^j\Gamma(\alpha
j+1)}{\Gamma(\alpha j-h)\Gamma(\nu j+1)}s_{k,h}\log^{\alpha
j-h-1}\left(1+A(1-u)\right)\right\}
\end{align*}
because $\frac{\alpha j-h}{\Gamma(\alpha
j-h+1)}=\frac{1}{\Gamma(\alpha j-h)}$ and $s_{k,0}=0$; then we get
\begin{align*}
\frac{d^{k+1}}{du^{k+1}}&E_{\nu,1}\left(-\log^\alpha\left(1+A(1-u)\right)\right)\\
=&\frac{(-A)^{k+1}}{(1+A(1-u))^{k+1}}\left\{-k\sum_{j\geq
0}\sum_{h=1}^k\frac{(-1)^j\Gamma(\alpha j+1)}{\Gamma(\alpha
j-h+1)\Gamma(\nu j+1)}s_{k,h}\log^{\alpha
j-h}\left(1+A(1-u)\right)\right.\\
&+\left.\sum_{j\geq 0}\sum_{m=1}^{k+1}\frac{(-1)^j\Gamma(\alpha
j+1)}{\Gamma(\alpha j-m+1)\Gamma(\nu j+1)}s_{k,m-1}\log^{\alpha
j-(m-1)-1}\left(1+A(1-u)\right)\right\}\\
=&\frac{(-A)^{k+1}}{(1+A(1-u))^{k+1}}\sum_{j\geq
0}\left\{\sum_{h=1}^k\frac{(-1)^j\Gamma(\alpha j+1)}{\Gamma(\alpha
j-h+1)\Gamma(\nu j+1)}(-ks_{k,h}+s_{k,h-1})\right.\\
&\left.\cdot\log^{\alpha
j-h}\left(1+A(1-u)\right)+\frac{(-1)^j\Gamma(\alpha
j+1)}{\Gamma(\alpha j-(k+1)+1)\Gamma(\nu j+1)}s_{k,k}\log^{\alpha
j-(k+1)}\left(1+A(1-u)\right)\right\},
\end{align*}
and \eqref{eq:derivative-by-induction} holds for $k+1$ because
$-ks_{k,h}+s_{k,h-1}=s_{k+1,h}$ and $s_{k,k}=s_{k+1,k+1}=1$.
$\Box$\\

In view of the next Proposition
\ref{prop:spacetimefractionalNB-eq-for-pgf} some preliminaries are
needed. Firstly let $(O)_\alpha$ be the operator defined by
\begin{equation}\label{eq:def-generalization-jaradetal}
(O)_\alpha f(z):=\left\{\begin{array}{ll}
\frac{1}{\Gamma(n-\alpha)}\int_{\frac{1-a}{b}}^z\log^{n-1-\alpha}\left(\frac{a+bz}{a+b\tau}\right)\left[\left(\left(\frac{a}{b}+\tau\right)
\frac{d}{d\tau}\right)^nf(\tau)\right]\frac{b}{a+b\tau}d\tau&\ \mbox{if}\ \alpha\in(n-1,n)\\
\left(\left(\frac{a}{b}+z\right)\frac{d}{dz}\right)^nf(z)&\
\mbox{if}\ \alpha=n
\end{array}\right.
\end{equation}
where $z>\frac{1-a}{b}$. Here, for the moment, we are assuming
that $\alpha>0$ and $n$ is an integer value. Thus, for
$\alpha\in(n-1,n)$, this operator can be formally considered as
the regularized Caputo-like fractional power of the operator
$\left(\frac{a}{b}+z\right)\frac{d}{dz}$. Indeed it can be found
from the definition of Caputo fractional derivative of order
$\alpha$, by means of the simple transformation
$z\mapsto\log(\frac{a}{b}+z)$. Moreover we observe that, if $a=0$
and $b=1$, \eqref{eq:def-generalization-jaradetal} coincides with
the Caputo-like regularized Hadamard fractional derivative
recently introduced in \cite{JaradAbdeljawadBaleanu}.

In what follows we focalize our attention on the case
$\alpha\in(0,1)$ and, in view of the proof of Proposition
\ref{prop:spacetimefractionalNB-eq-for-pgf}, we check that
\begin{equation}\label{eq:roberto2}
(O)_\alpha E_{\alpha,1}(-\gamma\log^\alpha(a+bz))=-\gamma
E_{\alpha,1}(-\gamma\log^\alpha(a+bz)).
\end{equation}
In fact, by \eqref{eq:def-generalization-jaradetal}, for
$\beta>-1$ we have
\begin{align*}
(O)_\alpha\log^\beta(a+bz)=&\frac{1}{\Gamma(1-\alpha)}
\int_{\frac{1-a}{b}}^z\log^{-\alpha}\left(\frac{a+bz}{a+b\tau}\right)\left[\left(\left(\frac{a}{b}+\tau\right)\frac{d}{d\tau}\right)
\log^\beta(a+b\tau)\right]\frac{b}{a+b\tau}d\tau\\
=&\frac{1}{\Gamma(1-\alpha)}
\int_{\frac{1-a}{b}}^z(\log(a+bz)-\log(a+b\tau))^{-\alpha}\cdot\beta\cdot\frac{\log^{\beta-1}(a+b\tau)}{a+b\tau}bd\tau
\end{align*}
and, after some computations with the change of variable
$y=\frac{\log(a+b\tau)}{\log(a+bz)}$, we obtain
$$(O)_\alpha\log^\beta(a+bz)=\frac{\beta}{\Gamma(1-\alpha)}\log^{\beta-\alpha}(a+bz)\int_0^1(1-y)^{-\alpha}y^{\beta-1}dy,$$
and therefore
\begin{equation}\label{eq:roberto1}
(O)_\alpha\log^\beta(a+bz)=\frac{\Gamma(\beta+1)}{\Gamma(\beta+1-\alpha)}\log^{\beta-\alpha}(a+bz);
\end{equation}
then, by \eqref{eq:roberto1} and some computations, we get
$$(O)_\alpha E_{\alpha,1}(-\gamma\log^\alpha(a+bz))=\sum_{k=1}^\infty\frac{(-\gamma)^k\log^{\alpha k-\alpha}(a+bz)}{\Gamma(\alpha k-\alpha+1)}
=-\gamma\sum_{k=0}^\infty\frac{(-\gamma)^k\log^{\alpha
k}(a+bz)}{\Gamma(\alpha k+1)},$$ which meets \eqref{eq:roberto2}.

\begin{proposition}\label{prop:spacetimefractionalNB-eq-for-pgf}
Assume that $r=1$ and let $\{G_{N_\rho^{\nu,\nu}(t)}:t\in [0,T]\}$
be the probability generating functions in Proposition
\ref{prop:spacetimefractionalNB-pgf} with $\alpha=\nu$. Then we
have the following results.\\
(i) (Case $\rho=1$) Let $(O)_{\nu,1}$ be the operator in
\eqref{eq:def-generalization-jaradetal} with $a=\frac{1}{p}$ and
$b=\frac{p-1}{p}$; then
$$(O)_{\nu,1}G_{N_1^{\nu,\nu}(t)}(u)=-G_{N_1^{\nu,\nu}(t)}(u)+1-\frac{\frac{1}{q(t)}-1}{\frac{1}{p}-1}\ (\mbox{for all}\ u\in(1,1/(1-p))).$$
(ii) (Case $\rho=0$) Let $(O)_{\nu,0}$ be the operator in
\eqref{eq:def-generalization-jaradetal} with $a=\frac{1}{q(t)}$
and $b=\frac{q(t)-1}{q(t)}$; then
$$(O)_{\nu,0}G_{N_0^{\nu,\nu}(t)}(u)=-G_{N_0^{\nu,\nu}(t)}(u)\ (\mbox{for all}\ u\in(1,1/(1-q(t)))).$$
(iii) In both cases (i) and (ii) we have
$G_{N_\rho^{\nu,\nu}(t)}(\frac{1-a}{b})=1$.
\end{proposition}
\noindent\emph{Proof.} We start with (i). For $\alpha=\nu\in(0,1)$
we have
\begin{align*}
(O)_{\nu,1}G_{N_1^{\nu,\nu}(t)}(u)=&\frac{\frac{1}{q(t)}-1}{\frac{1}{p}-1}(O)_{\nu,1}E_{\nu,1}\left(-\log^\nu\left(\frac{1-(1-p)u}{p}\right)\right)\\
=&-\frac{\frac{1}{q(t)}-1}{\frac{1}{p}-1}E_{\nu,1}\left(-\log^\nu\left(\frac{1-(1-p)u}{p}\right)\right)
=-G_{N_1^{\nu,\nu}(t)}(u)+1-\frac{\frac{1}{q(t)}-1}{\frac{1}{p}-1}
\end{align*}
where (for $\rho=1$, $a=\frac{1}{p}$, $b=\frac{p-1}{p}$ and
$\gamma=1$) we have used Proposition
\ref{prop:spacetimefractionalNB-pgf} and, for the second equality,
\eqref{eq:roberto2}. Note that we have $u\in(1,1/(1-p))$ because
$G_{N_1^{\nu,\nu}(t)}(u)$ is finite for $|u|<\frac{1}{1-p}$ (see
Proposition \ref{prop:spacetimefractionalNB-pgf} with $\rho=1$)
and $\frac{1-a}{b}=1$. For $\alpha=\nu=1$ it is easy to check with
some computations that
$$\left(\frac{1}{p-1}+u\right)\frac{d}{du}G_{N_1^{1,1}(t)}(u)=-G_{N_1^{1,1}(t)}(u)+1-\frac{\frac{1}{q(t)}-1}{\frac{1}{p}-1}.$$
by Proposition \ref{prop:spacetimefractionalNB-pgf} (in fact we
have $\frac{a}{b}=\frac{1}{p-1}$).\\
We proceed similarly for (ii). For $\alpha=\nu\in(0,1)$ we have
\begin{align*}
(O)_{\nu,0}G_{N_1^{\nu,\nu}(t)}(u)=&(O)_{\nu,0}E_{\nu,1}\left(-\log^\nu\left(\frac{1-(1-q(t))u}{q(t)}\right)\right)\\
=&-E_{\nu,1}\left(-\log^\nu\left(\frac{1-(1-q(t))u}{q(t)}\right)\right)
=-G_{N_0^{\nu,\nu}(t)}(u)
\end{align*}
where (for $\rho=0$, $a=\frac{1}{q(t)}$, $b=\frac{q(t)-1}{q(t)}$
and $\gamma=1$) we have used Proposition
\ref{prop:spacetimefractionalNB-pgf} and, for the second equality,
\eqref{eq:roberto2}. Note that we have $u\in(1,1/(1-p))$ arguing
as we did for the proof of (i). For $\alpha=\nu=1$ it is easy to
check with some computations that
$$\left(\frac{1}{q(t)-1}+u\right)\frac{d}{du}G_{N_0^{1,1}(t)}(u)=-G_{N_0^{1,1}(t)}(u)$$
by Proposition \ref{prop:spacetimefractionalNB-pgf} (in fact we
have $\frac{a}{b}=\frac{1}{q(t)-1}$).\\
Finally (iii) trivially holds because we always have
$G_{N_\rho^{\alpha,\nu}(t)}(1)=1$ (even if $\alpha\neq\nu$) and,
in both cases (i) and (ii), $\frac{1-a}{b}=1$. $\Box$

\section{On weighted processes}\label{sec:weighted-processes}
In this section we consider $\{N_\rho^w(t):t\in [0,T]\}$ where
$$N_\rho^w(t):=\sum_{n=1}^{M_g^w}1_{[0,t]}(X_n^{F,\rho})$$
and the probability mass function of the random variable $M_g^w$
is given by
\begin{equation}\label{eq:weighted-definition}
P(M_g^w=k)=\frac{P(M_g=k)w(k)}{\mathbb{E}[w(M_g)]}\ (\mbox{for
all}\ k\geq 0)
\end{equation}
for some nonnegative numbers (weights) $\{w(k):k\geq 0\}$ such
that
$$\mathbb{E}[w(M_g)]:=\sum_{r=0}^\infty w(r)P(M_g=r)\in (0,\infty);$$
then we are referring to the concept of weighted probability mass
function (see e.g. \cite{JohnsonKotzKemp}, p. 90, and the
references cited therein).

We remark that $M_g^w$ has the same distribution of $M_g$ if
$w(k)=1$ (for all $k\geq 0$). More in general we have the
following well-known property of the weighted probability mass
functions: if we consider \lq\lq proportional weights\rq\rq
$$\{w(k):k\geq 0\}\propto\{\tilde{w}(k):k\geq 0\},$$
i.e. if, for some $c>0$, we have $w(k)=c\tilde{w}(k)$ (for all
$k\geq 0$), then we have the same weighted probability mass
function.

The aim of this section is to illustrate the \lq\lq weighted
version structure\rq\rq\ for the probability mass function of
$N_\rho^w(t)$ for each $t\in(0,T]$, i.e.
\begin{equation}\label{eq:weights-depending-on-t-definition}
P(N_\rho^w(t)=k)=\frac{P(N_\rho(t)=k)w(k,t)}{\mathbb{E}[w(N_\rho(t),t)]}\
(\mbox{for all}\ k\geq 0)
\end{equation}
for some weights $\{w(k,t):k\geq 0\}$ which depend on $t\in(0,T]$
(obviously we have $w(k,T)=w(k)$ for all $k\geq 0$, i.e.
\eqref{eq:weights-depending-on-t-definition} meets
\eqref{eq:weighted-definition} when $t=T$). Moreover we give the
corrected version of some formulas in
\cite{BorgesRodriguesBalakrishnan}.

\begin{proposition}\label{prop:weighted-choiceofweights}
We set
\begin{align*}
q(k|n,F(t),\rho)&:=(1-\rho)\bincoeff{n}{k}F^k(t)(1-F(t))^{n-k}\\
&+\rho F^{k/n}(t)(1-F(t))^{1-k/n}1_{\{0,n\}}(k)\ (\mbox{for all}\
k\in\{0,1,\ldots,n\}).
\end{align*}
Then, for all $t\in(0,T]$, we have
\begin{equation}\label{eq:weighted-choiceofweights}
w(k,t)\propto\frac{\sum_{n=k}^\infty q(k|n,F(t),\rho)P(M_g=n)w(n)}
{\sum_{n=k}^\infty q(k|n,F(t),\rho)P(M_g=n)}\ (\mbox{for all}\
k\geq 0).
\end{equation}
\end{proposition}
\noindent\emph{Proof.} By (7) in
\cite{BorgesRodriguesBalakrishnan} we have the following
generalization of \eqref{eq:basic-pmf-rho=0}:
\begin{equation}\label{eq:basic-pmf-rho}
P(N_\rho(t)=k)=\sum_{n=k}^\infty q(k|n,F(t),\rho)P(M_g=n)\
(\mbox{for all}\ k\geq 0).
\end{equation}
Moreover, by \eqref{eq:basic-pmf-rho} (with $N_\rho^w(t)$ and
$M_g^w$ in place of $N_\rho(t)$ and $M_g$) and
\eqref{eq:weighted-definition}, we obtain
$$P(N_\rho^w(t)=k)=\frac{\sum_{n=k}^\infty q(k|n,F(t),\rho)P(M_g=n)w(n)}{\mathbb{E}[w(M_g)]}.$$
Then \eqref{eq:weights-depending-on-t-definition} and the last
equality yield
\begin{align*}
w(k,t)=&\frac{\mathbb{E}[w(N_\rho(t),t)]}{P(N_\rho(t)=k)}\cdot\frac{\sum_{n=k}^\infty q(k|n,F(t),\rho)P(M_g=n)w(n)}{\mathbb{E}[w(M_g)]}\\
\propto&\frac{\sum_{n=k}^\infty
q(k|n,F(t),\rho)P(M_g=n)w(n)}{P(N_\rho(t)=k)}.
\end{align*}
We conclude the proof by taking into account
\eqref{eq:basic-pmf-rho} for the denominator in the last
expression. $\Box$\\

Now the correction of (17) and (18) in
\cite{BorgesRodriguesBalakrishnan}:
$$\mbox{Cov}(N_\rho(t),N_\rho(s))=\lambda s\{1+\lambda\rho(1-t)\}$$
and
$$\mbox{Cov}(N_\rho(t)-N_\rho(s),N_\rho(s))=-\lambda^2\rho s(t-s).$$
We also present the corrected version of the displayed formula in
Example 4.1 in \cite{BorgesRodriguesBalakrishnan}. We refer to
\eqref{eq:basic-pmf} in this note and, in order to have a strict
connection with the presentation in
\cite{BorgesRodriguesBalakrishnan}, we consider $t\in[0,1]$ in
place of $F(t)$ with $t\in[0,T]$. We have to choose
$$P(N_0(t)=k)=\frac{(\lambda t)^k}{k!}e^{-\lambda t}\left(1-t+\frac{k}{\lambda}\right)\ (\mbox{for
all}\ k\geq 0)$$ for the case $\rho=0$ (see a displayed formula in
Section 3.1 in \cite{BalakrishnanKozubowski}) and
$$P(M_g=k)=\left\{\begin{array}{ll}
\frac{\lambda^{k-1}}{(k-1)!}e^{-\lambda}&\ \mbox{if}\ k\geq 1\\
0&\ \mbox{if}\ k=0;
\end{array}\right.$$
then we get
\begin{align*}
P(N_\rho(t)=k)&=(1-\rho)\frac{(\lambda t)^k}{k!}e^{-\lambda t}\left(1-t+\frac{k}{\lambda}\right)+
\rho\left\{(1-t)1_{k=0}+t\cdot\frac{\lambda^{k-1}}{(k-1)!}e^{-\lambda}\cdot 1_{k\geq 1}\right\}\\
&=\left\{\begin{array}{ll}
(1-\rho)e^{-\lambda t}(1-t)+\rho(1-t)&\ \mbox{if}\ k=0\\
(1-\rho)\frac{(\lambda t)^k}{k!}e^{-\lambda
t}\left(1-t+\frac{k}{\lambda}\right)+\rho
t\frac{\lambda^{k-1}}{(k-1)!}e^{-\lambda}&\ \mbox{if}\ k\geq 1,
\end{array}\right.
\end{align*}
which is the corrected version of the displayed formula in Example
4.1 in \cite{BorgesRodriguesBalakrishnan}.

\paragraph{Acknowledgements.} We thank the referee for some useful
comments and Federico Polito for Figure \ref{fig}.

\end{document}